\documentclass[a4paper,11pt,english,leqno]{amsart}
\usepackage{amsmath, amsthm, amsfonts,amssymb, mathrsfs}
\usepackage{verbatim, graphicx, ifthen, enumitem}
\usepackage[T1]{fontenc}
\usepackage[dvipsnames]{xcolor}
\usepackage{hyperref}
\usepackage{cleveref}
\usepackage{tikz}
\usepackage{bm}
\usepackage{caption}
\usepackage{subcaption}
\usepackage[english]{babel}
\usetikzlibrary{calc}
\usepackage[normalem]{ulem}

\newcommand\myshade{85}
\colorlet{mylinkcolor}{violet}
\colorlet{mycitecolor}{YellowOrange}
\colorlet{myurlcolor}{Aquamarine}

\hypersetup{
  linkcolor  = mylinkcolor!\myshade!black,
  citecolor  = mycitecolor!\myshade!black,
  urlcolor   = myurlcolor!\myshade!black,
  colorlinks = true,
}

\newtheorem{thm}{Theorem}[subsection]

\newtheorem{lem}[thm]{Lemma}
\newtheorem{prop}[thm]{Proposition}
\theoremstyle{definition}
\newtheorem{defn}[thm]{Definition}

\theoremstyle{remark}
\newtheorem{rem}[thm]{Remark}

\numberwithin{equation}{subsection}
\numberwithin{figure}{subsection}

\newcommand{\diff}{\mathrm{d}}
\newcommand{\C}{{\mathbb C}}
\newcommand{\R}{{\mathbb R}}

\newcommand{\Z}{{\mathbb Z}}

\newcommand{\dA}{{\mathrm{dA}}}

\newcommand{\imag}{\mathrm{i}}
\newcommand{\e}{\mathrm{e}}

\newcommand{\Stwohalf}{\sph^2_{\frac12}}

\newcommand{\Tope}{{\mathbf T}}
\newcommand{\ordo}{\mathrm{o}}
\newcommand{\Ordo}{\mathrm{O}}

\newcommand\sph{\mathbb{S}}
\newcommand\surf{\mathscr{S}}
\newcommand{\balpha}{\boldsymbol\alpha}

\def\E#1{\mathcal{E}_{#1}}

\DeclareMathOperator{\re}{Re}
\DeclareMathOperator{\im}{Im}

\makeatletter
\let\@@citation@@=\citation

\renewcommand{\citation}[1]{\@@citation@@{#1}%
\@for\@tempa:=#1\do{\@ifundefined{cit@\@tempa}%
  {\global\@namedef{cit@\@tempa}{}}{}}%
}

\def\@lbibitem[#1]#2#3\par{%
  \@ifundefined{cit@#2}{}{\item[\@biblabel{#1}\hfill]}%
  \if@filesw
      {\let\protect\noexpand
       \immediate
       \write\@auxout{\string\bibcite{#2}{#1}}}\fi\ignorespaces
  \@ifundefined{cit@#2}{}{#3}}
\def\@bibitem#1#2\par{%
  \@ifundefined{cit@#1}{}{\item}%
  \if@filesw \immediate\write\@auxout
    {\string\bibcite{#1}{\the\value{\@listctr}}}\fi\ignorespaces
  \@ifundefined{cit@#1}{}{#2}}
\makeatother

\begin{document}

%
\title{A Bombieri-type inequality and equidistribution of points}

\author{Ujué Etayo}
\address{Etayo: Departamento de Matemáticas, CUNEF Universidad, Leonardo Prieto Castro, 2.
Ciudad Universitaria, 28040 Madrid, Spain}
\email{ujue.etayo@cunef.edu}

\author{Haakan Hedenmalm}
\address{
Hedenmalm: Department of Mathematics\\
KTH Royal Institute of Technology\\
S--10044 Stockholm\\
Sweden
\\
Department of Mathematics and Statistics\\
University of Reading\\
RG6 6AX Reading, England\\
Department of Mathematics and Computer Science
\\
St Petersburg State University
\\
St Petersburg, Russia}

\email{haakan00@gmail.com}

\author{Joaquim Ortega-Cerd\`a}
\address{
Ortega-Cerd\`a: Departament de Matem\`atiques i Inform\`atica \\
Universitat de  Barcelona \& Centre de Recerca Matemàtica \\
Gran Via 585, 08007, Barcelona \\
 Spain}
\email{jortega@ub.edu}

 
\thanks{The research of Uju\'e Etayo was supported by the Austrian Science Fund FWF project F5503
(part of the Special Research Program (SFB) Quasi-Monte Carlo Methods: Theory and Applications),
by the Grant PID2020-113887GB-I00 funded by
MCIN/AEI/10.13039/501100011033
and by the
starting grant from BBVA associated to the prize Jos\'e Luis Rubio de Francia. The research
of Haakan Hedenmalm was supported by
Vetenskapsr\aa{}det (VR grant 2020-03733), the Leverhulme trust (grant VP1-2020-007), and 
by grant 075-15-2021-602 of the Government of the Russian Federation for the state support
of scientific research, carried out under the supervision of leading scientists.
Joaquim Ortega-{Cerd\`a} was supported by grants PID2021-123405NB-I00
and CEX2020-001084-M by the Agencia Estatal de Investigaci\'on and by the Departament de
Recerca i Universitats, grant 2021 SGR 00087 and ICREA Academia}
 
\begin{abstract} 
In recent work, Etayo introduces a new Bombieri-type inequality for monic polynomials.
Here we reinterpret this new 
inequality as
{a more general} integral inequality involving
the Green function for the sphere.
This rather geometric interpretation allows for generalizations of the basic
inequality, involving fractional zeros while also opening up the possibility to extend
the setting to general compact Riemann surfaces.
We derive a sharp form of these generalized Bombieri-type inequalities for the case of
the sphere and the torus. These inequalities involve a quantity we call the
\emph{packing number}, which in turn is inspired by the geometric zero packing problems
considered by Hedenmalm in the context of the asymptotic variance of the
Bergman projection of a bounded function. 
As for the torus, we introduce analogs of polynomials (pseudopolynomials)
based on the classical Weierstrass $\sigma$ function, and we explain how such
pseudopolynomials fit in with the extended geometric Bombieri-type inequality.
The sharpness of the packing number bound on the torus involves the construction of a
lattice configuration on the torus for any given integer number of points.
The corresponding bound for the sphere instead relies on the existence of
well-conditioned polynomials in the sense of Shub and Smale.
\end{abstract}

\maketitle

\section{Introduction and background material}\label{sec_INT}

\subsection{Basic notation}
\label{subsec-1.1}
We write $\R$ for the real line and $\C$ for the complex plane. Moreover, we write
$\mathbb{S}=\C\cup\{\infty\}$ for the extended complex plane (the Riemann sphere). 
For a complex variable $z=x+\imag y\in\C$, let 
\[
\dA(z):=\diff x\wedge\diff y
\]
denote the usual area measure
and 
\begin{equation*}
\dA_{\mathbb{S}}(z)=\frac{\dA(z)}{(1+|z|^2)^2}
\end{equation*}
the area measure on the Riemann sphere. 

\subsection{Stereographic projection on the sphere }
Let $\Stwohalf$ be the sphere in $\R^3$ of radius $\frac12$ centered at
$\left(0,0,\frac12\right)$, with the top point $(0,0,1)$ as the north pole. 
The points of the sphere $\Stwohalf$ and of the extended complex plane
$\sph=\C\cup\{\infty\}$ are identified using stereographic
projection, except that the north pole is identified with the added point
at infinity. The stereographic projection
maps $\pi_\sph:\Stwohalf\to\sph=\C\cup\{\infty\}$,
\begin{equation}\label{eq:stepro}
\pi_\sph:\,(a,b,c)\mapsto\frac{a+\imag b}{1-c},
\end{equation}
where $a^2+b^2+\left(c-\frac12\right)^2=\frac14$ and, in the reverse direction,
\[
\pi_\sph^{-1}:\,z=x+\imag y\mapsto \left(\frac{x}{1+|z|^2},\frac{y}{1+|z|^2},
\frac{|z|^2}{1+|z|^2}\right).
\] 
The Euclidean distance in $\R^3$ between two points $q,q'\in\Stwohalf\subset\R^3$
with $z,z'\in\C$ and
$z=\pi_\sph(q)$, $z'=\pi_\sph(q')$, is the quantity
\begin{equation}\label{eq:sphmetric}
|q-q'|:=\diff_\sph(z,z')=
\frac{|z-z'|}{(1+|z|^2)^{\frac12}(1+|z'|^2)^{\frac12}}.
\end{equation}
Next, letting $z'$ approach $z$, we obtain from \cref{eq:sphmetric} the spherical
length element 
\[
|\diff q|=\frac{|\diff z|}{1+|z|^2},
\]
which is in agreement with our previous definition of the spherical area element
$\dA_\sph$.
It is worth observing that the distance between the points $q,q'$ in
\eqref{eq:sphmetric} is at most
$1$, and that distance $1$ is attained precisely when the points are
\emph{antipodal}: $q'=-q$. In
terms of the complex coordinates $z,z'$ in \eqref{eq:sphmetric}, this means that
$\bar zz'=-1$, and we then say that $z,z'$ are \emph{antipodal in the complex plane}.

Using a homothetic transformation of ratio $2$ and with center at the north pole
we may obtain instead the sphere centered at the origin of radius $1$, and when needed, we
denote this sphere by $\sph_{1}^{2}$.
Observe that the stereographic projection defined in equation \eqref{eq:stepro}
also maps the sphere
$\sph_{1}^{2}$ onto the extended complex plane.
\bigskip

\subsection{The Bombieri-Weyl norm}\label{sec_Bombi}

Let $N\ge1$ be an integer, and consider the Bombieri-Weyl space 
$\mathrm{BW}_N$ of bivariate homogeneous polynomials of degree $N$, that is,
the $\C$-linear vector space of dimension $N+1$ consisting of all polynomials of
the form
\begin{equation}\label{eq:qpol1}
Q(z,w)=\sum_{j=0}^{N} a_j\, z^{j}w^{N-j},
\end{equation}
where the coefficients $a_j$ are complex numbers. 
The Weyl norm of $q$ (also known as Bombieri norm or Kostlan norm) is given by
\begin{equation}\label{eq:bombi}
\|Q\|_N^2:=\sum_{j=0}^{N}\binom{N}{j}^{-1}\,|a_j|^2,
\end{equation}
which makes clear the dependence on $N$. 
Alternatively, the Weyl norm may be expressed in the form
\[
\|Q\|_N^2=\frac{N+1}{\pi}
\int_{\mathbb{P}\C^{1}}\frac{|Q(\xi)|^2}{|\xi|^{2N}}\diff\sigma(\xi),
\]
where the integration is made with respect to the volume form $\sigma$ arising from
the standard Riemannian structure of the projective space $\mathbb{P}\C^{1}$. 
The norm can also be expressed as
\[
\|Q\|_N^2=(N+1)\int_{\partial B} |Q(\xi)|^2 d\mu(\xi),
\]
where the integral is over the unit sphere $\partial B$ in $\C^2$ and the measure $\mu$ is the
spherical measure normalized so that $\mu(\partial B) = 1$.
Some classical references on the Weyl norm are \cite{BBEM,10.2307/2690843}.
The homogeneous polynomials of degree $N$ of the form \cref{eq:qpol1} are in a
one-to-one correspondence
with the univariate polynomials 
\begin{equation*}
Q(z):=Q(z,1)=\sum_{j=0}^{N} a_j\, z^{j},
\end{equation*}
the recovery formula being
\begin{equation*}
Q(z,w)=w^NQ(z/w).
\end{equation*}
In terms of the univariate polynomial, the Weyl norm may be written as
\begin{equation}\label{eq_02}
\|Q\|_N^2
=\sum_{j=0}^{N}\binom{N}{j}^{-1}\,|a_j|^2
=
\frac{N+1}{\pi}
\int_{\C}\frac{|Q(z)|^2}{(1+|z|^{2})^N}\dA_\sph(z).
\end{equation}
For instance, if $Q(z)=(z-a)^N$, we calculate that
\begin{equation}
\|Q\|_N^2=\|(z-a)^N\|_N^2=\sum_{j=0}^N\binom{N}{j}^{-1}\binom{N}{j}^2|a|^{2(N-j)}=(1+|a|^2)^N.
\label{eq:powernorm}  
\end{equation}

\subsection{The Bombieri inequality for bivariate polynomials}

If $Q \in \mathrm{BW}_N$ factors as $Q=Q_1Q_2$, where $Q_j$ has degree $N_j$ for $j=1,2$,
with $N=N_1+N_2$, the Bombieri inequality asserts that
\begin{equation}\label{eq:Bomb0}
\frac{N_1!N_2!}{N!}\,\|Q_1\|^2_{N_1}\|Q_2\|^2_{N_2}\le\|Q_1Q_2\|^2_{N_1+N_2}
=\|Q\|^2_N.
\end{equation}
On the other hand, it is known that the Weyl norm is submultiplicative
(see \cite{Beauzamy90}), i.e.,
\begin{equation}\label{algebranorm}
\|Q\|^2_N \le\|Q_1\|^2_{N_1}\|Q_2\|^2_{N_2}.
\end{equation}
In fact, the upper bound expressed in \eqref{algebranorm} is sharp. We just
choose $Q(z)=(z-a)^N$ and
and $Q_j(z) = (z-a)^{N_j}$ for $j=1,2$ and any point $a\in\mathbb{C}$, and
observe that, according to
\eqref{eq:powernorm}, we have
\begin{equation*}
\|Q\|^2_N =(1+|a|^2)^N=(1+|a|^2)^{N_1}(1+|a|^2)^{N_2}= \|Q_1\|^2_{N_1}\|Q_2\|^2_{N_2}.
\end{equation*} 
As for lower bound \eqref{eq:Bomb0} expressed by the Bombieri inequality, it is
sharp too, in the following sense:
given two natural numbers $N_1,N_2$ then there exist two polynomials $Q_1,Q_2$ of degree
respectively $N_1,N_2$ and such that
\begin{equation}
\frac{N_1!N_2!}{N!}\,\|Q_1\|^2_{N_1}\|Q_2\|^2_{N_2}
=
\|Q_1Q_2\|^2_{N_1+N_2}
=
\|Q\|^2_N,
\label{eq:sharpBomb0}
\end{equation}
see the original statement due to Beauzamy, Bombieri, Enflo and Montgomery in
\cite{BBEM}.
In our rather simple situation of bivariate homogeneous polynomials, which we
reduce to ordinary
univariate polynomials, it is easy to exhibit explicit examples of such polynomials
$Q_1,Q_2$ and $Q=Q_1Q_2$.
If we choose $Q_1(z)=z^{N_1}$ and $Q_2(z)=1$, then
\[
\|Q_1\|_{N_1}^2=\|Q_2\|^2_{N_2}=1,\quad \|Q\|_{N}^2=\|Q_1Q_2\|_N^2=
\binom{N}{N_1}^{-1}=\frac{N_1!N_2!}{N!},
\]
and hence equality is attained in \eqref{eq:Bomb0}. If we analyze this example carefully,
we can use M\"obius invariance on the sphere to see that \eqref{eq:sharpBomb0} holds more generally
for $Q=Q_1Q_2$, if $Q_1(z)=(z-a)^{N_1}$ and $Q(z)=(z-b)^{N_2}$, provided that the points
$a,b\in\C$ are antipodal in the plane:
$\bar a b=-1$. So, at the intuitive level, sharpness in the Bombieri estimate
\eqref{eq:Bomb0} corresponds to the zeros of $Q_1$ and $Q_2$ each being concentrated
to a point, but the two points repel each other, while sharpness in the
submultiplicative estimate \eqref{eq:powernorm} instead only involves concentration
(attraction).

Naturally, we may iterate the Bombieri inequality \eqref{eq:Bomb0}, and obtain, if
$Q=Q_1\cdots Q_k$ holds while $Q_j$ is of degree $N_j$, and $N=N_1+\ldots+N_k$,
\begin{multline}\label{eq:iterBomb1}
\frac{N_1!\cdots N_k!}{N!}\,\|Q_1\|^2_{N_1}\cdots\|Q_k\|^2_{N_k}
\le\|Q_1\cdots Q_k\|^2_{N_1+\cdots+N_k}=\|Q\|^2_N
\\
\le\|Q_1\|^2_{N_1}\cdots
\|Q_k\|^2_{N_k}.
\end{multline}
From the intuition developed in the study of the sharpness of the
Bombieri inequality \eqref{eq:Bomb0},
we would expect the that in each step when we factor $Q$ the sharpness would
occur from a different
configuration of zeros, so that e.g., the optimal configuration of zeros
for the factorization
$Q=Q_1Q_2Q_3$ should be that each polynomial $Q_j$ has concentrated zeros at a
single point, but that the
three (multiple) zeros should all repel each other. However, this configuration
does not agree with what
should be {optimal} at the earlier step, when, e.g., $Q=Q_1(Q_2Q_3)$. This means
that we cannot expect
the iterated inequality \eqref{eq:iterBomb1} to be sharp as soon as it involves
three or more nontrivial factors.
Let us bring this to extreme, and factor $Q$ down to first degree polynomials,
as permitted by the fundamental theorem of algebra.
We then have $N_j=1$ for all $j=1,\ldots,N$, and the assertion of
\eqref{eq:iterBomb1} -- combined with the
iterated submultiplicativity based on \eqref{eq:powernorm} -- reduces to
\begin{equation}\label{eq:iterBomb2}
\frac{1}{N!}\|Q_1\|^2_{1}\cdots\|Q_N\|^2_{1}
\le\|Q\|^2_N\le\|Q_1\|^2_{1}\cdots\|Q_N\|^2_{1}.
\end{equation}
The loss of sharpness in the left-hand side inequality of \eqref{eq:iterBomb2}
for three or more factors was apparently first noted by Carlos Beltr\'an and
Uju\'e Etayo.
In fact, in \cite{Etayo}, Etayo established the {asymptotically} considerably
better inequality
\begin{equation}\label{eq:iterBomb2.1}
c_N\|Q_1\|^2_{1}\cdots\|Q_N\|^2_{1}
\le\|Q\|^2_N\le\|Q_1\|^2_{1}\cdots\|Q_N\|^2_{1}
\end{equation}
where the parameter $c_N$ takes the form
\begin{equation}
\label{quan:cN}
c_N:=\frac{N+1}{d_{N} \e^{N}}
\end{equation}
and where $d_{N}$ satisfies $0 < \gamma_0 < d_{N} < 1$ for a certain constant $\gamma_0$,
the smallest such constant for which the inequality holds.
The quantity $d_{N}$ is conjectured to be of the form
\begin{equation*}
d_N = K + o(1)\quad\text{as}\,\,\,N\to+\infty,
\end{equation*}
with $K$ a universal constant, see Problem 5.1 of \cite{Etayo} and also \cite{BF}.
The constant $c_N$ appearing in the inequality \eqref{eq:iterBomb2.1} is considerably
better than the corresponding constant $1/N!$ in the inequality \eqref{eq:iterBomb2}
at least asymptotically as $N\to+\infty$. Moreover, while not known with high precision,
the constant $c_N$ is  sharp in the following sense.
According to Theorem 4.5 of \cite{Etayo}, we have that for each $N=1,2,3,\ldots$
there exist
some first degree polynomials $Q_1,\cdots ,Q_N$ such that with $Q=Q_1\cdots Q_N$,
\begin{equation*}
c_N \|Q_1\|^2_{1}\cdots\|Q_N\|^2_{1}
=
\|Q\|^2_N.
\end{equation*}


%



\section{Main results and structure of the paper}\label{sec_MR}


\subsection{Integral form of inequality \eqref{eq:iterBomb2.1}}

We will now view the iterated Bombieri-type inequality \eqref{eq:iterBomb2.1} from the
prism of the auxiliary function
\begin{equation}\label{eq:Usphere}
U_{\sph}(z,w)=\frac12+\frac{1}{2}\log\frac{|z-w|^2}{(1+|z|^2)(1+|w|^2)}, \quad z,w\in\C,
\end{equation}
which has a natural interpretation in terms of potential theory on the sphere. 
For points $w_1,\ldots,w_N \in \C$, we consider the associated polynomial with
those zeros $Q(z) = \prod_{j=1}^{N}Q_j(z)$, where $Q_j(z)=z-w_j$ are the first-degree
factors.
We observe that by the integral definition of the Weyl norm \eqref{eq_02}, we have that
\begin{equation}
\frac{1}{\pi}
\int_{\C} \exp\left( 2U_{\sph}(z,w_j)\right) \dA_{\mathbb{S}}(z)
=
\frac{\e}{2}\,\frac{\|Q_j\|_1^2}{1+|w_j|^2}=\frac{\e}{2} 
\label{eq:normaffine}
\end{equation}
for $j=1,\ldots,N$, and, more generally,
\begin{multline*}
\frac{1}{\pi}\int_{\C}\exp\big(2 U_\sph(z,w_1)+\cdots+2 U_\sph(z,w_N)\big) \dA_\sph(z)
=
\frac{\e^N \|Q\|_{N}^{2}}{(N+1)\prod_{j=1}^{N}(1+|w_j|^2)}.
\end{multline*}
From this point of view, the inequality \eqref{eq:iterBomb2.1} asserts that
\begin{multline*}
\frac{\e^N c_N\|Q_1\|^2_{1}\cdots\|Q_N\|^2_{1}}{(N+1)\prod_{j=1}^{N}(1+|w_j|^2)}
\leq
\frac{\e^N \|Q \|_{N}^{2}}{(N+1)\prod_{j=1}^{N}(1+|w_j|^2)}
\\
\leq
\frac{\e^N \|Q_1\|^2_{1}\cdots\|Q_N\|^2_{1}}{(N+1)\prod_{j=1}^{N}(1+|w_j|^2)},
\end{multline*}
which is the same as
\begin{multline}\label{eq:05}
c_N\prod_{j=1}^{N}
\frac{2}{\pi}
\int_{\C} \exp\left(2 U_{\sph}(z,w_j)\right) \dA_{\mathbb{S}}(z)
\\
\leq
\frac{N+1}{\pi}\int_{\C}\exp\big(2 U_\sph(z,w_1)+\cdots+2 U_\sph(z,w_N)\big) \dA_\sph(z)
\\
\leq
 \prod_{j=1}^{N}
\frac{2}{\pi}
\int_{\C} \exp\left(2 U_{\sph}(z,w_j)\right) \dA_{\mathbb{S}}(z).
\end{multline}
Moreover, in view of \eqref{eq:normaffine},
the expressions on the left-hand and right-hand sides of \eqref{eq:05} are known,
so that \eqref{eq:05} reduces to
\begin{multline}\label{eq:05.1}
\frac{c_N\e^N}{N+1}
\leq
\frac{1}{\pi}\int_{\C}\exp\big(2 U_\sph(z,w_1)+\cdots+2 U_\sph(z,w_N)\big) \dA_\sph(z)
\\
\leq
\frac{\e^N}{N+1}.
\end{multline}
The optimality of these two inequalities means that 
\begin{multline}
\inf_{w,\ldots,w_N}
\frac{1}{\pi}\int_{\C}\exp\big(2 U_\sph(z,w_1)+\cdots+2 U_\sph(z,w_N)\big) \dA_\sph(z)  
\\
=\frac{c_N\e^N}{N+1}=\frac{1}{d_N},
\label{eq:05.2}
\end{multline}
where the parameter $d_N$ is given by \eqref{quan:cN}, and that
\begin{equation}
\label{eq:05.3}
\sup_{w_1,\ldots,w_N}
\frac{1}{\pi}\int_{\C}\exp\big(2 U_\sph(z,w_1)+\cdots+2 U_\sph(z,w_N)\big) \dA_\sph(z)
=\frac{\e^N}{N+1}.
\end{equation}
Optimality in the latter inequality \eqref{eq:05.3} is achieved when the points
all coalesce: $w_1=\ldots=w_N$. {On} the other hand, optimality in
the first inequality \eqref{eq:05.2} arises when the points are in some sense optimally
dispersed over the sphere. 
Inspired by the first part of the latter inequality, we define the
\emph{packing number}
\begin{equation}
\varTheta_{N,\beta}(\sph):=
\inf_{w_1,\ldots,w_N}
\frac{1}{\pi}\int_{\C}\exp\big(\beta U_\sph(z,w_1)+\cdots+\beta U_\sph(z,w_N)\big)
\dA_\sph(z),
\label{eq:varthetasph}
\end{equation}
for real $\beta>-2$. The quantity $\beta$ amounts to a charge parameter, with equal
``electrostatic'' charges at all the points $w_1,\ldots,w_N$. 
%
Our first result consists in obtaining relevant bounds on this quantity
$\varTheta_{N,\beta}(\sph)$, independent of parameter $N$.

\begin{thm}\label{mt:1}
Let the packing number  $\varTheta_{N,\beta}(\sph)$ be defined as above.
Then $\varTheta_{N,0}(\sph)=1$, while for $\beta>0$, we have the estimates
\begin{equation*}
1
\le
\varTheta_{N,\beta}(\sph)
\le 
K^{\beta},
\end{equation*}
where $K>1$ is a constant.
\end{thm}

The proof of this theorem is supplied in Subsection \ref{**}.

\begin{rem}
By H\"older's inequality, the quantity
\[
\big(\varTheta_{N,\beta}(\sph)\big)^{1/\beta}
\]
is an increasing function of $\beta$ in the interval $0<\beta<+\infty$.
Consequently, optimizing over $N=1,2,3,\ldots$, we find that
\[
\sup_N \big(\varTheta_{N,\beta}(\sph)\big)^{1/\beta}
\]
is increasing in $\beta$ as well. The optimal constant $K$ in the theorem is
hence given by 
\[
K=\sup_{N,\beta} \big(\varTheta_{N,\beta}(\sph)\big)^{1/\beta}=\lim_{\beta\to+\infty}
\sup_N \big(\varTheta_{N,\beta}(\sph)\big)^{1/\beta}.  
\]
\end{rem}


\subsection{Generalization to compact Riemann surfaces}

The function introduced in the preceding subsection $U_\sph(z,w)$
can be understood as a version of the Green function for the 
sphere $\sph$. Of course, it is well-known that on a compact Riemann surface $\surf$
there is no function which is harmonic except for a single logarithmic ``pole''. 
However, there exists instead a \emph{bipolar Green function} $L_\surf(z,w,w')$, with two
logarithmic poles $w,w'$ of opposite sign, for $w\ne w'$. 
If in the  corresponding charts the function has the local expansions
$L_\surf(z,w,w')=\log|z-w|+\Ordo(1)$  for $z\sim w$ and $L(z,w,w')=-\log|z-w'|+\Ordo(1)$ for
$z\sim w'$, it is then uniquely determined up to an additive constant. 
To fix the constant  as well, we will require that 
\[
\int_{\surf}L_\surf(z,w,w')\,\dA_\surf(z)=0,
\]  
where $\dA_\surf$ is the canonical area measure on $\surf$,  having constant Gaussian
curvature. 
There is also a \emph{unipolar Green function} $U_{\surf}(z,w)$, which we obtain from
$L_\surf(z,w,w')$ by simply averaging over $w'$ with respect to the canonical area measure. 
In other words, if $w'\in \surf$ is thought of as random with density proportional to the
surface area measure $\dA_\surf$, then the expected value of  $L_\surf(z,w,w')$ equals the
function $U_\surf(z,w)$. It is worth noticing that the bipolar Green function can be recovered
from the unipolar Green function via the simple formula
\[
L_\surf(z,w,w')=U_\surf(z,w)-U_\surf(z,w').
\]
It is easy to see that the function  $U_\surf(z,w)$ when $\surf=\sph$ then coincides with
the expression $U_\sph(z,w)$ given by \cref{eq:Usphere}. 

\begin{prop}\label{mt4}
Let $\surf$ be a compact Riemann surface and let $U_\surf(z,w)$ be its corresponding unipolar
Green function. Let $\dA_\surf(z)$ be the area form associated to the metric in $\surf$
and let $|\surf|_A>0$ denote its area with respect to the canonical area mesure $\dA_\surf$.
Then, for any $\alpha >0$, there is a constant $C>1$ such that for any probability
measure $\rho$ on $\surf$ we have that
$$
1 \le  
\frac{1}{|\surf|_{\mathrm{A}}}\int_{\surf}
\exp\left(\alpha\int_{\surf}U_\surf(z,w)d\rho(w)\right)
\dA_\surf(z) \le C^\alpha.
$$
\end{prop}

The proof of this proposition is supplied in Subsection \ref{ss:simple}.

In the general case of a compact Riemann surface, it is not clear how to bound uniformly
the infimum of
\[
\frac{1}{|\surf|_{\mathrm{A}}}
\int_{\surf} \exp\left(\alpha N\int_{\surf} U_\surf(z,w)d\rho_N(w)\right)\dA_\surf(z)
\]
taken over all the discrete probability measures $\rho_N$ which place point masses of
equal mass at $N$ points of $\surf$.
However, we are able to do so when the surface $\surf$ is a torus.

\subsection{The case of the torus}
Any complex torus $\surf$ is modelled by $\surf=\C/\Lambda$, where $\Lambda$ is a
discrete lattice in $\C$. Without loss of generality, we assume that the 
covolume is $\frac12\pi$ (this is the area of any fundamental rhombus in $\C/\Lambda$).
The unipolar
Green function $U_{\C/\Lambda}(z,w)$ is then given by the explicit formula
\begin{equation}
U_{\C/\Lambda}(z,w)=\log|\sigma_\ast(z-w)|-|z-w|^2-A_\Lambda, \quad z,w\in\C,
\label{eq:unipolar1.1}
\end{equation}
where $\sigma$ is the usual Weierstrass sigma function while $\sigma_*$ is the modified
Weierstrass sigma function, as expressed by the relation \eqref{eq:sigma-star} below.
Moreover, the constant $A_\Lambda$ is the quantity 
\[
A_\Lambda:=\frac{2}{\pi}\int_{\C/\Lambda}(\log|\sigma_\ast(z)|-|z|^2)\,\dA(z).
\]
In a fashion analogous to that of the sphere, given points
$w_1,\ldots,w_N \in \C$ and $\beta>-2$, we define the associated packing
number
\begin{multline}\label{eq:defVtorus}
\varTheta_{N,\beta}(\C/\Lambda):=
\\
\inf_{w_1,\ldots,w_N}
\frac{2}{\pi}
\int_{\C/\Lambda}\exp\big(\beta U_{\C/\Lambda}(z,w_1)+\cdots+\beta U_{\C/\Lambda}
(z,w_N)\big)\dA(z).
\end{multline}
We have the following theorem.

\begin{thm}\label{mt:2}
Let $\varTheta_{N,\beta}(\C/\Lambda)$ be defined as above and suppose that
$\beta >0$. Then
\begin{equation*}
1
\le
\varTheta_{N,\beta}(\C/\Lambda)
\le 
K^{\beta}
\end{equation*}
holds, where the constant $K>1$ only depends on the lattice $\Lambda$.
\end{thm}


\subsection{Pseudopolynomials on the torus}
\label{ss:pseudopols}

In the case of the sphere $\sph$, we have an equivalence between integrals involving
the unipolar Green function $U_{\mathbb{S}}(z,w)$ and Bombieri norms of polynomials. As
for the case of the torus $\C/\Lambda$, it is natural to consider the
\emph{first degree monic pseudopolynomial on the torus $\C/\Lambda$}
\[
\pi_\alpha(z):=\e^{\imag2\im \bar\alpha z}
\,\e^{-|z-\alpha|^2}\sigma_\ast(z-\alpha),
\qquad \alpha,z \in \C,
\]
which has a root at the point $\alpha$ modulo the lattice in $\C/\Lambda$.
It follows from the definition of the unipolar Green function \eqref{eq:unipolar1.1}
that the following identity holds:
\begin{equation}
|\pi_\alpha(z)|^\beta=\e^{\beta A_\Lambda}\exp(\beta U_{\C/\Lambda}(z,\alpha)).
\label{eq:pibeta1}
\end{equation}
The \emph{monic pseudo\-polynomials of degree $N$ on the torus $\C/\Lambda$} 
are declared to be the functions 
\[
\Pi_{\balpha}(z):=\pi_{\alpha_1}(z)\cdots\pi_{\alpha_N}(z),
\]
where $\balpha=(\alpha_1,\ldots,\alpha_N)$ are the \emph{roots} of $\Pi_{\balpha}$. 
For $\beta>0$, we let $\|\cdot \|_{L^\beta}$ denote the $L^{\beta}$-norm in
$\C/\Lambda$ with $\dA$ as the underlying measure (which by assumption puts total mass
$\frac12\pi$ on $\C/\Lambda$), and find that 
\begin{multline*}
\left|\left|\pi_{\alpha}\right|\right|_{L^\beta}^{\beta}
=\int_{\C/\Lambda}|\pi_\alpha(z)|^\beta\dA(z)=
e^{\beta A_{\Lambda}}
\int_{\C/\Lambda}\exp(\beta U_{\C/\Lambda}(z,\alpha))\,\dA(z)
\\
=\frac{\pi}{2}\,\e^{\beta A_\Lambda}E(\beta,\Lambda),
\end{multline*}
where we use translation invariance and introduce the notation
\begin{equation}
E(\beta,\Lambda):=\frac{2}{\pi}\int_{\C/\Lambda}\exp(\beta U_{\C/\Lambda}(z,0))\,\dA(z).
\label{eq:EbetaN}
\end{equation}
Similarly, the $L^\beta$-norm of the monic pseudopolynomial $\Pi_{\balpha}$ is
\begin{multline}
\label{eq:Lbetanorm2.01}
\left|\left|\Pi_{\balpha}\right|\right|_{L^\beta}^{\beta}
=\int_{\C/\Lambda}|\Pi_{\balpha}(z)|^\beta\dA(z)
\\
=\e^{N\beta A_{\Lambda}}
\int_{\C/\Lambda}\exp\bigg(\beta \sum_{j=1}^{N}U_{\C/\Lambda}(z,\alpha_j)\bigg)\,\dA(z),
\end{multline}
which can be estimated from above and below.


\begin{prop}\label{mt:3}
Let $\alpha_{1},\ldots,\alpha_{N} \in \C$ and $\beta>0$. Then we have the following estimate
of the $L^\beta$-norm of the monic pseudopolynomial $\Pi_{\balpha}=\prod_j \pi_{\alpha_j}$:  
\begin{equation*}
\frac{\pi}{2}\,\e^{N\beta A_{\Lambda}}\,\varTheta_{N,\beta}(\C/\Lambda) 
\le
\|\Pi_{\balpha}\|_{L^\beta}^{\beta}
\le 
\frac{\pi}{2}\,\e^{N\beta A_{\Lambda}}\,{E(\beta N,\Lambda)}.
\end{equation*} 
Here, the quantity $\varTheta_{N,\beta}(\C/\Lambda)$ is given by the
relation \eqref{eq:defVtorus}, while the expression $E(\beta N,\Lambda)$ is
as in \eqref{eq:EbetaN}.
\end{prop}


\subsection{Organization of the paper}

In Section \ref{sec_INT}, we supply the background on iterated Bombieri 
inequalities, and later, in Section \ref{sec_MR}, we present the framework of
more general geometric Bombieri-type inequalities, while outlining our main results.
Furthermore, in Section \ref{sec:IF}, we present the inequality \eqref{eq:iterBomb2.1}
as a particular case of a more general geometrically inspired inequality. In the same
section, we also obtain Theorem \ref{mt:1} and Proposition \ref{mt4}.
In \cref{sec_torus} we present the integral form of the Bombieri type inequality
for the torus and we write down the proof of Theorem \ref{mt:2}.
Finally, in Section \ref{sec_TP}, the main properties of the pseudopolynomials
for the torus are studied, and we supply the proof of {Proposition} \ref{mt:3}.



\section{An inequality on the sphere} \label{sec:IF}

\subsection{The unipolar Green function on the sphere}

For $z,w \in \C$, let $U_{\sph}$ be the function  
\begin{equation*}
U_{\sph}(z,w)=\frac12+\frac{1}{2}\log\frac{|z-w|^2}{(1+|z|^2)(1+|w|^2)}
\end{equation*}
which assumes values in the extended interval $\left[-\infty,\frac12\right]$ and
vanishes on average over the sphere, 
\begin{equation}
\label{eq_03}
\int_\C U_\sph(z,w)\dA_\sph(z)=0,\qquad w\in\C.
\end{equation}
A direct computation involving a rotation of the sphere $\sph$ reveals that
\begin{equation}\label{eq:elint1}
\frac{1}{\pi}\int_\C \exp(\alpha U_\sph(z,w))\dA_\sph(z)=
\frac{2\,\e^{\alpha/2}}{2+\alpha}
\end{equation}
independently of $w$ provided that the integral converges (which happens precisely
when $\alpha>-2$).
\begin{prop}\label{prop_mu_gen}
For a real parameter $\alpha>0$ and a Borel probability measure $\rho$ on $\C$,
we have that
\[
1 \le \frac{1}{\pi} 
\int_{\mathbb C} \exp\left(\alpha \int_{\mathbb C } U_{\sph}(z,w)\diff\rho(w)\right)
\dA_\sph(z) \le \frac{2\, \e^{\alpha/2}}{2+\alpha }.
\]
\end{prop}

\begin{proof}
It follows from Jensen's inequality that
\begin{multline}
\frac{1}{\pi} 
\int_{\mathbb C} \exp\left(\alpha \int_{\mathbb C }U_{\sph}(z,w)\diff\rho(w)\right) \dA_\sph(z)
\\
\le
\frac{1}{\pi} 
\int_{\mathbb C} \int_{\mathbb C } \exp(\alpha U_{\sph}(z,w))\diff\rho(w) \dA_\sph(z)
\\
=\frac{1}{\pi} 
\int_{\mathbb C} \int_{\mathbb C } \exp(\alpha U_{\sph}(z,w))\dA_\sph(z)\,\diff\rho(w) =
\frac{2\, \e^{\alpha/2}}{2+\alpha }\int_{\mathbb C} \diff\rho{(w)}
 =\frac{2\, \e^{\alpha/2}}{2+\alpha },
\end{multline}
by Fubini's theorem and the identity \eqref{eq:elint1}.
As for the lower bound, we again use Jensen's inequality and Fubini's theorem in combination
with \eqref{eq_03} to obtain that
\begin{multline}\label{eq_unif}
\frac{1}{\pi} 
\int_{\mathbb C} \exp\left(\alpha\int_{\mathbb C} U_{\sph}(z,w)d\rho(w)\right) \dA_\sph(z) 
\\
\ge 
\exp\left(\frac{\alpha}{\pi}
\int_{\mathbb C} \int_{\mathbb C } U_{\sph}(z,w)d\rho(w) \dA_\sph(z)
\right) = \e^0 = 1.
\end{multline}
\end{proof}

\begin{rem}
In the special case $\diff\rho=\dA_\sph$, it holds that 
\[
\int_{\mathbb C } U_{\sph}(z,w)\diff\rho(w) = 0,
\]
and hence then the lower bound \eqref{eq_unif} is actually an equality.
If instead $\diff\rho$ happens to be a discrete probability measure which is similar to
$\dA_\sph$ in that it is as evenly spread over the sphere as possible, how much of this
phenomenon remains?
\end{rem}

\subsection{Integral form of the inequality 
\texorpdfstring{\eqref{eq:iterBomb2.1}}{(\ref{eq:iterBomb2.1})}}\label{IFI}

Given $N$ points $w_1,\ldots,w_N \in\C$, we let
\[
\rho_N = \frac 1{N}\sum_{j= 1}^N  \delta_{w_j}
\]
be the associated probability measure on $\mathbb C$ placing equal mass $1/N$ at
each of the points.
It is immediate from Proposition \ref{prop_mu_gen} (with $\alpha=N\beta$) that
\begin{equation}\label{eq:iterBomb3}
1 \le \frac{1}{\pi} 
\int_{\mathbb C} \exp\left(\beta N\int_{\mathbb C } U_{\sph}(z,w)d\rho_N(w)\right)
\dA_\sph(z) \le \frac{2\, \e^{\beta N/2}}{2+\beta N}.
\end{equation}
The integral expression in \eqref{eq:iterBomb3} is maximized when all the
points $w_1,\ldots,w_N$ coincide, in which case we may apply \eqref{eq:elint1}:
\[
\max_{w_1,\ldots,w_N}
\frac{1}{\pi}\int_{\C}\exp\big(\beta U_\sph(z,w_1)+\cdots+\beta U_\sph(z,w_N)\big)
\dA_\sph(z)=\frac{2\,\e^{\beta N/2}}{2+N\beta}.
\]
In particular, for $\beta = 2$, we recover the right-hand side inequality
of \eqref{eq:iterBomb2.1}. 
As for the left-hand side inequality of \cref{eq:iterBomb2.1}, it is instead
optimal when the points $w_1,\ldots,w_N$ are well spread out over the sphere. 
We recall the quantity $\varTheta_{N,\beta}(\sph)$ from the preceding section given
by \eqref{eq:varthetasph},
\begin{equation*}
\varTheta_{N,\beta}(\sph)=\inf_{w_1,\ldots,w_N}
\frac{1}{\pi}\int_{\C}\exp\big(\beta U_\sph(z,w_1)+\cdots+\beta U_\sph(z,w_N)\big)
\dA_\sph(z).
\end{equation*}
In view of Proposition \ref{prop_mu_gen} we have the estimate from below
\begin{equation}
\varTheta_{N,\beta}(\sph)\ge
1.
\label{eq:lowerbound01}
\end{equation}
Moreover, H\"older's inequality gives that for $0<\beta<\beta'$,
\begin{multline*}
\frac{1}{\pi}\int_{\C}\exp\big(\beta U_\sph(z,w_1)+\cdots+\beta U_\sph(z,w_N)\big)
\dA_\sph(z)
\\
\le \Bigg(
\frac{1}{\pi}\int_{\C}\exp\big(\beta' U_\sph(z,w_1)+\cdots+\beta' U_\sph(z,w_N)\big)
\dA_\sph(z)\Bigg)^{\beta/\beta'},
\end{multline*}
and consequently,
\[
\varTheta_{N,\beta}(\sph)^{1/\beta}\le\varTheta_{N,\beta'}(\sph)^{1/\beta'},
\qquad 0<\beta<\beta'.
\]
By the improved H\"older's inequality of \cite{Hed2}, we get to compare
$\beta$ and $\beta'=2\beta$:
\[
\varTheta_{N,\beta}(\sph)^2\le(1-\varrho_{N,\beta}(\sph))\varTheta_{N,2\beta}(\sph),
\]
for $\beta>0$. Here, the constant 
\begin{multline*}
\varrho_{N,\beta}(\sph):=
\\
\inf
\,\,\frac{1}{\pi}\int_{\C} \left( 
\frac{a\,|z - w_1|^\beta \cdots |z - w_N|^\beta}{(1 + |z|^2)^{\frac{N\beta}{2}}
(1 + |w_1|^2)^{\frac{\beta}{2}} \cdots (1 + |w_N|^2)^{\frac{\beta}{2}} } - 1
\right)^2
\dA_\sph(z) 
> 0
\end{multline*}
is associated with so-called spherical zero packing, as defined in \cite{Hed2}.
The infimum runs over all $a\in\R$ and all point configurations
$w_1,\ldots,w_N\in\C$.
By iteration, it follows that
\begin{equation*}
\varTheta_{N,\beta}(\sph)
\ge\prod_{j=1}^{+\infty}(1-\rho_{N,2^{-j}\beta}(\sph))^{-2^{j-1}},
\end{equation*}
where each factor in the right-hand side product is $>1$.
In \cite{Hed2} it is conjectured that 
\begin{equation*}
\lim\limits_{N \longrightarrow +\infty}
\varrho_{N,\beta} (\sph)
\longrightarrow \varrho_{\beta}(\mathbb{C}),
\end{equation*}
where $\varrho_{\beta}(\mathbb{C})$ is a monotonically increasing function
with
\begin{equation*}
\lim\limits_{\beta \longrightarrow 0^{+}}
\varrho_{\beta}(\mathbb{C})
 =
0
\qquad
\text{and}
\qquad
\lim\limits_{\beta \longrightarrow +\infty}
\varrho_{\beta}(\mathbb{C})
 =
 1.
\end{equation*}
It follows from \cite{Etayo} that
$\varTheta_{N,\beta}$ is bounded uniformly in $N$, for fixed $\beta=2$. 
In the next section, we extend this result to all $\beta >0$.

\subsection{The uniform bound on the packing number
$\varTheta_{N,\beta}(\sph)$
}
\label{**}
In \cite{BEMOC}, it is shown that there exists a real constant $\gamma>0$
such that for all $N=1,2,3,\ldots$, we can find a collection of points
$X_{N} = \{ x_1,\cdots ,x_{N} \}\in \Stwohalf$
with the properties that
for each $q \in \Stwohalf$
we have that $\mathrm{dist}(q,X_{N}) \leq \gamma/\sqrt{N}$,
where $\mathrm{dist}(q,X_{N})$ denotes the Euclidean distance from $q$ to
the set $X_{N}$ in $\R^3$ while at the same time
\begin{equation}\label{BEMOC1}
B_1\,\frac{\sqrt{N}\,\mathrm{dist}(q,X_{N})}{\e^{N/2}}
\leq  \prod_{j=1}^{N} |q-x_{j}|
\leq B_2\,
\frac{\sqrt{N}\,\mathrm{dist}(q,X_{N})}{\e^{N/2}},
\end{equation}
for two positive constants $B_1<B_2$ independent of $N$. 
The results as formulated in Theorem 1.11 of \cite{BEMOC} are expressed for the
sphere $\sph_1^2$ of radius $1$ but all the distances in $\Stwohalf$ are just
are just half the distances in the $\sph_1^2$ model and thus easily accounted for. 
Let us explore what this entails as regards the upper bound on the quantity
$\varTheta_{N,\beta}(\sph)$. 
According to the definition \eqref{eq:varthetasph}, we have that
\begin{multline*}
\varTheta_{N,\beta}(\sph)
=
\inf_{w_1,\ldots,w_N}
\frac{1}{\pi}\int_{\C}\exp\big(\beta U_\sph(z,w_1)+\cdots+
\beta U_\sph(z,w_N)\big)\,\dA_\sph(z)
\\
=
\inf_{w_1,\ldots,w_N}
\frac{\e^{\beta N/2}}{\pi}
\int_{\C}
\prod_{j=1}^N\frac{|z-w_j|^\beta}
{(1+|z|^2)^{\beta/2}(1+|w_j|^2)^{\beta/2}}
\,\dA_\sph(z)
\\
=
\inf_{p_1,\ldots,p_N}
\frac{\e^{\beta N/2}}{\pi}
\int_{\Stwohalf}
\prod_{j=1}^N|q - p_{j}|^{\beta}
\,
\diff\sigma(q),
\end{multline*}
where $d\sigma(q)$ stands for the Lebesgue area measure on
the sphere $\Stwohalf$ (which has total mass $\pi$) and
$p_1,\ldots,p_{N}\in \Stwohalf$ are the preimages
of the points $w_1,\ldots,w_{N}$ by the stereographic projection,
i.e. $p_{j} = \pi_{\sph}^{-1}(w_{j})$ for $j=1,\ldots, N$.
Then, by the existence of the set of points $X_N=\{x_1,\ldots,x_N\}$
in $\Stwohalf$ with \eqref{BEMOC1} we find that
\begin{multline}
\varTheta_{N,\beta}(\sph)=
\inf_{p_1,\ldots,p_N}
\frac{\e^{\beta N/2}}{\pi}
\int_{\Stwohalf}
\prod_{j=1}^{N}
|q - p_{j}|^{\beta}
\diff\sigma(q)
\\
\leq
\frac{\e^{\beta N/2}}{\pi}
\int_{\Stwohalf}
\prod_{j=1}^{N} |q-x_{j}|^{\beta}
d\sigma(q)
\leq
\frac{\e^{\beta N/2}N^{\beta/2}B_2^\beta}{\pi}
\int_{\sph_{1}^{2}}
\frac{(\mathrm{dist}(q,X_{N}))^\beta}{\e^{\beta N/2}}
d\sigma(q)
\\
=
\frac{N^{\beta/2}B_2^{\beta}}{\pi}
\int_{\Stwohalf} (\mathrm{dist}(q,X_{N}))^\beta
d\sigma(q)
\leq
\gamma^\beta B_2^\beta.
\label{eq:optconditioning}
\end{multline}

\begin{proof}[Proof of Theorem \ref{mt:1}]
The theorem follows with $K=\gamma B_2$ by \eqref{eq:optconditioning}.
\end{proof}

\subsection{Spherical packing for $\beta = 2$ and Shub-Smale's
condition number.}
In the particular case $\beta = 2$, we have a separate
characterization of the quantity $\varTheta_{N,2}(\sph)$.
As a first step, we observe that if $P_N$ is the monic degree
$N$ polynomial with zeros at $w_1,\ldots,w_N$, that is,
\[
P_N(z)=\prod_{j=1}^{N}(z-w_j),
\]
then
\begin{equation}
\varTheta_{N,2}(\sph)=\inf_{w_1,\ldots,w_N}
\frac{\|P_N\|_N^2}{N+1}\prod_{j=1}^{N}(1+|w_j|^2)^{-1},  
\label{eq:packing1.01}
\end{equation}
in view of the definition \eqref{eq:varthetasph} and the
integral expression for the Hilbert space norm $\|\cdot\|_N$.
While this is a nice formula, we will try to find another
which involves more natural geometrically motivated expressions.
We have $N$ points $X_N=\{x_1,\ldots,x_N\}$ on the sphere
$\Stwohalf$, and the points $w_j=\pi_\sph(x_j)\in\C\cup\{\infty\}$
arise by stereographic projection.
The logarithmic energy of the set of points
$X_{N} = \{x_{1}, \ldots  ,x_{N} \}$ is given by 
\[
\E{\log}(X_{N})= \sum_{j,k:j \neq k} \log\frac{1}{|x_{j}-x_{k}|}
=\frac12\sum_{j,k:j \neq k}
\log\frac{(1+|w_j|^2)(1+|w_k|^2)}{|w_{j}-w_{k}|^2},
\]
where the summation is over $j,k$ in $\{1,\ldots,N\}$.
The logarithmic energy is related to the transfinite diameter and the
capacity of the given set in the sense of classical potential theory,
see the books \cite{doohovskoy2011foundations} or \cite{Saff}. 
We also need the concept of \emph{condition number}, as conceived by
Shub and Smale \cite{SS93}.
Given a polynomial $p$ of degree $\le N$ and one of its roots,
call it $\zeta\in\C$, the condition number $\mu_N(p,\zeta)$ of $p$
at the root $\zeta$ measures how quickly the root $\zeta$ moves as we
perturb the polynomial $p$. According to \cite{SS93}, the formula on
p. 6 in combination with Lemma 1, the condition number --
associated with the polynomial $p$ of degree $\le N$ and its root
$\zeta\in\C$ -- is given by
\begin{equation}
\mu_N(p,\zeta)=\sqrt{N}\|p\|_N\,\frac{(1+|\zeta|^2)^{(N-2)/2}}{|p'(\zeta)|}.
\label{eq:munumber1}
\end{equation}
So, at multiple roots, the condition number would be infinite.
We notice first that by \eqref{eq:munumber1}, 
\begin{equation*}
\mu_N(P_N,w_j)=\sqrt{N}\|P_N\|_N\,\frac{(1+|w_j|^2)^{(N-2)/2}}{|P_N'(w_j)|}.
\end{equation*}
so that
\begin{equation}
\label{eq:prodcondition2}
\prod_{j=1}^N\mu_N(P_N,w_j)=
\frac{{N}^{N/2}\|P_N\|_N^N}{|P_N'(w_1)\cdots P_N'(w_N)|}
\prod_{j=1}^N(1+|w_j|^2)^{(N-2)/2}.
\end{equation}
Next, since
\begin{equation*}
P_N'(w_j)=\prod_{k:k\ne j}(w_j-w_k),
\end{equation*}
it follows that
\begin{equation*}
\prod_{j=1}^{N}|P_N'(w_j)|=\prod_{j,k:k\ne j}|w_j-w_k|,
\end{equation*}
and, consequently,
\begin{multline}
\exp\big(-\E{\log}(X_{N})\big)= \prod_{j,k:j \neq k}
\frac{|w_{j}-w_{k}|}{(1+|w_j|^2)^{\frac12}(1+|w_k|^2)^{\frac12}}
\\
=\prod_{j=1}^{N}|P_N'(w_j)|\,\prod_{j=1}^N(1+|w_j|^2)^{1-N}.
\label{eq:logenergy2}
\end{multline}
It now follows from a combination of \eqref{eq:prodcondition2}
and \eqref{eq:logenergy2} that
\begin{equation*}
\exp\big(-\E{\log}(X_{N})\big)\prod_{j=1}^N\mu_N(P_N,w_j) 
=N^{N/2}\|P_N\|_N^N\,\prod_{j=1}^N(1+|w_j|^2)^{-N/2},
\label{eq:logenergy2*}
\end{equation*}
so that
\begin{multline}
\exp\bigg(-\frac{2}{N}\,\E{\log}(X_{N})\bigg)
\prod_{j=1}^N(\mu_N(P_N,w_j))^{2/N} 
\\
=N\|P_N\|_N^2\,\prod_{j=1}^N(1+|w_j|^2)^{-1}.
\label{eq:logenergy3}
\end{multline}

\begin{prop}
Let $X_{N} = \{x_{1},\ldots,x_{N} \}$ be a collection of $N$ separate
points on the sphere $\Stwohalf$, and let $w_{j} =\pi_{\sph}(x_{j})$
be the corresponding complex points obtained by stereographic projection.
Moreover, let $P_N$ denote the monic polynomial with the complex points
as zeros, 
\[
P_{N}(z) = \prod_{j=1}^{N} (z-w_{j}).
\]  
The packing number for $N$ points and point charge $\beta=2$ is then
given by
\begin{equation*}
\varTheta_{N,2}(\sph)
=
\frac{1}{N(N+1)}
\inf_{w_1,\ldots,w_N}
\exp\bigg(-\frac{2}{N}\,\E{\log}(X_N)\bigg)
\prod_{j=1}^{N} \big(\mu_{N}(P_{N},w_{j})\big)^{2/N}.
\end{equation*}
\end{prop}

\begin{proof}
The formula is an immediate consequence of \eqref{eq:packing1.01}
and \eqref{eq:logenergy3}.
\end{proof}

\subsection{A simple integral bound for general compact
surfaces}
\label{ss:simple}
We obtain the following.

\begin{proof}[Proof of Proposition \ref{mt4}]
The proof follows along the lines of Proposition \ref{prop_mu_gen},
the only difference being that we need an upper bound for the integral 
\[
\int_{\surf}\exp\bigl( \alpha U_{\surf}(z,w)\bigr)\dA_\surf(z).
\]  
It is, however, well known (see \cite{Aubin}) that the unipolar Green function is
bounded from above, that is, $U_{\surf}(z,w)\le M$ holds for some positive constant
$M$, and hence
\[
\frac{1}{|\surf|_{\mathrm{A}}}\int_{\surf}\exp\bigl( \alpha U_{\surf}(z,w)\bigr)
\dA_\surf(z)\le \e^{\alpha M}=C^\alpha,
\]
where $C=\e^M>1$. 
\end{proof}

\section{A Bombieri-type inequality for the torus}
\label{sec_torus}

\subsection{Modelling the torus}

We model the torus as $\C/\Lambda$, where $\Lambda$ is a discrete lattice
(additive subgroup of $\C\cong\R^2$) having two $\R$-linearly
independent generators. So, $\Lambda$ is of the form  
\[
\Lambda=\omega_1\Z+\omega_2\Z
\]
in $\C$, where the generators are $\omega_1,\omega_2\in\C$ which are
$\R$-linearly independent.
After an appropriate rotation of the plane, we may suppose that 
$\omega_1\in\R$ with $\omega_1>0$, and that $\im \omega_2>0$.

\subsection{The Weierstrass functions $\wp$, $\zeta$,
and $\sigma$}
Associated with the given lattice $\Lambda$ there is the Weierstrass elliptic
function $\wp(z)$ which is $\Lambda$-periodic with a double pole at the
origin on the torus $\C/\Lambda$, with asymptotics
\[
\wp(z)=\frac{1}{z^2}+\sum_{\omega\in\Lambda\setminus\{0\}}
\bigg(\frac{1}{(z-\omega)^2}-\frac{1}{\omega^2}\bigg)=\frac{1}{z^2}+\Ordo(1),
\quad\text{ as}\,\,\,\,z\to0.
\]
We note that since $\omega\in\Lambda$ holds if and only if
$-\omega\in\Lambda$, the Weierstrass function $\wp$ is even: $\wp(-z)=\wp(z)$.
In connection with the Weierstrass $\wp$-function, we have the auxiliary
Weierstrass functions $\zeta$ and $\sigma$. The $\zeta$ function is, as a matter
of definition, the primitive of $-\wp$, and hence it has a simple pole with
residue $1$ at each point of $\Lambda$. It is then determined uniquely by the
additional requirement to be odd: $\zeta(-z)=-\zeta(z)$.
However, the periodicity of $\wp$ does not carry over to the primitive $\zeta$,
and instead we get integration constants $\lambda_1,\lambda_2$ such that for
$j=1,2$,
\begin{equation*}
\zeta(z+\omega_j)=\zeta(z)+\lambda_j,\qquad z\in \C.
\end{equation*}
Moreover, in view of the odd symmetry, we read off that
\[
\lambda_j=2\zeta\left(\tfrac{1}{2}\omega_j\right),\qquad j=1,2.
\]  
The $\sigma$ function is entire and its logarithmic derivative equals the
function $\zeta$.
The explicit formula is (see, e. g., \cite{Ahlfors}):
\[
\sigma(z):=z\prod_{0\ne\omega\in\Lambda}
\bigg(1-\frac{z}{\omega}\bigg)\exp\bigg\{
\frac{z}{\omega}+\frac{z^2}{2\omega^2}
\bigg\}.
\]  
As is clear from the product formula, the function $\sigma$ is entire with
simple zeros along $\Lambda$, 
and enjoys the 
periodicity-type formul\ae{}
\begin{equation}
\sigma(z+\omega_j)=-\sigma(z)
\exp\big(\lambda_j\left( z+\tfrac12{\omega_j}\right) \big),
\qquad j=1,2.
\label{eq-period1}
\end{equation}
We define the modified Weierstrass $\sigma$ function to be the function
$\sigma_\ast$ given by  
\begin{equation}
\sigma_\ast(z):=\e^{a z+b z^2}\sigma(z),
\label{eq:sigma-star}
\end{equation}
where $a,b\in\C$ are parameters which we need to determine. 
The property which we want to be fulfilled is that the associated function
\begin{equation}\label{eq-sigmafunction1}
\e^{-|z|^2}|\sigma_\ast(z)|=\e^{-|z|^2+\re(a z+b z^2)}|\sigma(z)|
\end{equation}
be $\Lambda$-periodic. 
This is possible precisely when the fundamental rhombus
$\mathcal{D}$ in $\C/\Lambda$ has area equal to $\frac1{2}\pi$. In terms
of the periods $\omega_1,\omega_2$, the area requirement amounts to 
\begin{equation}
\label{eq:areacond1}
\omega_1\im \omega_2=\frac{\pi}{2},
\end{equation}
which in turn connects with the classical Legendre relation 
$\lambda_1\omega_2-\lambda_2\omega_1=\imag2\pi$.
Under the area condition \eqref{eq:areacond1}, it is indeed possible
to specify values for the constants
$a$ and $b$ such that the function in \eqref{eq-sigmafunction1}
gets to be doubly periodic. We select
\begin{equation*}
a=0,\qquad b=1-\frac{\lambda_1}{2\omega_1}=
\frac{\bar\omega_2}{\omega_2}-\frac{\lambda_2}{2\omega_2}, 
\end{equation*}
and calculate that
\begin{equation}
\sigma_\ast(z)=-\e^{-2\bar\omega_j z-|\omega_j|^2}\sigma_\ast(z+\omega_j)
=-\Tope_{\omega_j}\sigma_\ast(z),\qquad
j=1,2,
\label{eq:sigmaper1}
\end{equation}
where $\Tope_\alpha$ denotes the \emph{Weyl translate}
\begin{equation}
\Tope_\alpha f(z):=\e^{-2\bar\alpha z-|\alpha|^2}f(z+\alpha),
\label{eq:focktran1}
\end{equation}
which acts unitarily on the corresponding Fock space:
\[
\int_\C|\Tope_\alpha f(z)|^{2}\,\e^{-2|z|^2}\dA(z)=
\int_\C|f(z)|^{2}\,\e^{-2|z|^2}\dA(z).
\]
The standard commutation relation for the Weyl translates reads
\begin{equation*}
\Tope_{\alpha'}\Tope_\alpha f(z)=\e^{-\imag2\im\bar\alpha\alpha'}
\Tope_{\alpha+\alpha'}f(z)
=\e^{-\imag4\im\bar\alpha\alpha'}\Tope_{\alpha}\Tope_{\alpha'} f(z).
\end{equation*}
It is elementary but a bit tedious to check that the
function $\e^{-|z|^2}|\sigma_\ast(z)|$ is well-defined and
continuous on the torus $\C/\Lambda$, and vanishes at
the origin, and we verify that the function
\[
U_{\C/\Lambda}(z,w)=U_{\C/\Lambda}(z-w,0)=\log|\sigma_\ast(z-w)|-|z-w|^2-A_\Lambda
\]
equals the unipolar Green function encountered earlier in the context
of the torus surface $\C/\Lambda$.
Here, $A_\Lambda$ is the number 
\[
A_\Lambda:=\frac{2}{\pi}\int_{\C/\Lambda}(\log|\sigma_\ast(z)|-|z|^2)\dA(z),
\]
which equilibrates things so that the average of the function
$U_{\C/\Lambda}(\cdot,w)$ on $\C/\Lambda$ vanishes: 
\[
\int_{\C/\Lambda}U_{\C/\Lambda}(z,w)\dA(z)=0.
\]
Let $w_{1},\ldots,w_{N}\in\C/\Lambda$ and $\beta >0$, then as a
corollary to Proposition \ref{mt4}, we obtain that
\begin{multline}
\label{eq:04}
1 \leq
\varTheta_{N,\beta}(\C/\Lambda)\leq
\\
\frac2{\pi}\int_{\C/\Lambda}
\exp\big(\beta U_{\C/\Lambda}(z,w_1)+\cdots+\beta U_{\C/\Lambda}(z,w_N)\big)
\,\dA(z)\le
E(\beta N,\Lambda).
\end{multline} 
where 
the integrals
\begin{equation*}
E(\beta,\Lambda):=\frac{2}{\pi}\int_{\C/\Lambda}
\exp(\beta U_{\C/\Lambda}(z,w))dA(z)
\end{equation*}
are independent of the point $w\in\C$, and
\[
\sup_{w_1,\ldots,w_N}
\frac{2}{\pi}\int_{\C/\Lambda}\exp\big(\beta U_{\C/\Lambda}(z,w_1)+\cdots+
\beta U_{\C/\Lambda}(z,w_N)\big)\, \dA(z)
=
E(\beta N,\Lambda)
\]
is attained when all the points $w_1,\ldots,w_N$ coincide. 


\subsection{The uniform bound on the packing number
for the torus}

We obtain the uniform boundedness of the packing number
$\varTheta_{N,\beta}(\C/\Lambda)$ independently of the number $N$  of point
charges $\beta>0$ by using the structure of the torus $\C/\Lambda$.
We construct a lattice $\Lambda_N$ which contains $\Lambda$ as a sublattice
with the property that $\Lambda_N/\Lambda$ consists of exactly $N$ points.  
Here, by "lattice" we mean a discrete additive subgroup of 
$\langle\C,+\rangle$. 
It {will} be an important feature that the typical distance 
between adjacent points in $\Lambda_N/\Lambda$ is of the order $\Ordo(N^{-1/2})$. 

Given $N=1,2,3,\ldots$, let $a=a_N=\lfloor\sqrt{N}\rfloor$ be the integer part
of $\sqrt{N}$. Moreover, let $b=b_N$ be the integer $b_N=N-a_N^2$, which falls in
the interval $0\le b_N<2\sqrt{N}$. We then define the
complex numbers $\omega_{j,N}$ for $j=1,2$ by the formula
\begin{equation}
  \begin{pmatrix}
\omega_{1,N}\\
\omega_{2,N}
\end{pmatrix}
:=\frac{1}{N}
\begin{pmatrix}
a_N &1 
\\
-b_N & a_N
\end{pmatrix}
\begin{pmatrix}
 \omega_{1}\\
\omega_{2} 
\end{pmatrix}
\label{eq:matrixrel1}
\end{equation}
where we recall that $\omega_1,\omega_2$ are the generators of the lattice
$\Lambda$, with $\omega_1>0$ and $\im\omega_2>0$. 
We introduce the lattice
\begin{equation}
\Lambda_N:=\omega_{1,N}\Z+\omega_{2,N}\Z,
\end{equation}
and since the inverse relationship to \eqref{eq:matrixrel1} reads
\begin{equation}
  \begin{pmatrix}
\omega_{1}\\
\omega_{2}
\end{pmatrix}
=
\begin{pmatrix}
a_N &-1 
\\
b_N & a_N
\end{pmatrix}
\begin{pmatrix}
 \omega_{1,N}\\
\omega_{2,N} 
\end{pmatrix}
\label{eq:matrixrel2}
\end{equation}
with integer entries in the matrix, it is immediate that
we have the lattice inclusion $\Lambda\subset\Lambda_N$.
It remains to count the number of points in $\Lambda_N/\Lambda$.
To this end, we may use that the area change is determined by
the determinant of the matrix in \eqref{eq:matrixrel2}, and if
we associate a small rhombic cell to each point of
$\Lambda_N/\Lambda$, we get that the total area of the torus $\C/\Lambda$
is given by the number of points in $\Lambda_N/\Lambda$ times
the area of $\C/\Lambda_N$. As the area change is also given by the
determinant, it follows that the total number of points in
$\Lambda_N/\Lambda$ equals
\[
\det\begin{pmatrix}
a_N &-1 
\\
b_N & a_N
\end{pmatrix}
=a_N^2+b_N=N.
\]


We have already discussed the unipolar Green function $U_\Lambda$ associated with
the lattice $\Lambda$. There is a corresponding unipolar Green function for the
lattice $U_{\Lambda_N}$. 

\begin{lem}\label{lem_01}
We list the $N$ points $\{\zeta_1,\ldots,\zeta_N\}$ in
$\Lambda_N/\Lambda\subset\C/\Lambda$.
Then we have the identity
\begin{equation*}
\sum_{j=1}^{N} 
U_{\C/\Lambda}(z,\zeta_j)
=
U_{\C/\Lambda_{N}}(z,0),
\qquad z\in\C/\Lambda.
\end{equation*}
\end{lem}

\begin{proof}
The function $\sum_{j=1}^{N} U_{\C/\Lambda}(\cdot,\zeta_j)$ is
seen to be $\Lambda_{N}$-periodic, and its Laplacian is
given by 
\begin{equation*}
\Delta 
\sum_{j=1}^{N} 
U_{\C/\Lambda}(\cdot,\zeta_j)
=
2\pi\sum_{j=1}^{N}
\delta_{\zeta_j} - 4N\quad\text{on}\,\,\,\,\C/\Lambda,
\end{equation*}
since we have that
\begin{equation*}
\Delta U_{\C/\Lambda}(\cdot,\zeta)
=
2\pi\,
\delta_{\zeta} - 4\quad\text{on}\,\,\,\,\C/\Lambda.
\end{equation*}
On the other hand, the function $U_{\C/\Lambda_{N}}(z,0)$ is also
$\Lambda_{N}$-periodic, and as such its Laplacian is given by 
\begin{equation*}
\Delta U_{\C/\Lambda_{N}}(z,0)
=
2\pi\sum_{\lambda \in \Lambda_{N}}
\delta_{\lambda} - 4N.
\end{equation*}
The increase in the background constant to $4N$ compared with $4$ for
$\C/\Lambda$ comes from the need to compensate for the smaller area of
$\C/\Lambda_N$. Comparing the two expressions on the torus, we find that
the two expressions must coincide up to an additive constant.
Moreover, since each of the two expressions vanish in the mean on $\C/\Lambda$,
the additive constant is in fact $0$ and hence the expressions coincide.
\end{proof}

\begin{lem}
Let $\Lambda'_N$ be the renormalized lattice
\[
\Lambda'_N :=\sqrt{N}\Lambda_N=\omega'_{1,N}\Z+\omega'_{2,N}\Z,
\]
where
\begin{equation}
  \begin{pmatrix}
\omega'_{1,N}\\
\omega'_{2,N}
\end{pmatrix}
:=\frac{1}{\sqrt{N}}
\begin{pmatrix}
a_N &1 
\\
-b_N & a_N
\end{pmatrix}
\begin{pmatrix}
 \omega_{1}\\
\omega_{2} 
\end{pmatrix}.
\label{eq:matrixrel1.5}
\end{equation}
If $\C/\Lambda'_N$ is the associated torus, then the associated
fundamental rhombus has area $\frac12\pi$, while the generating vectors
$\omega'_{j,N}$ are uniformly bounded in $N$. Moreover, the unipolar Green
function $U_{\C/\Lambda'_N}(\cdot,0)$ is uniformly bounded:
\[
U_{\C/\Lambda'_N}(z,0)\le M,\qquad z\in\C,
\]
holds for some constant $M>0$ independent of $N$.  
\label{lem:renormalized1}
\end{lem}

\begin{proof}
The lengths of the generating vectors are controlled,
\[
|\omega'_{1,N}|=|\omega_1|+\ordo(1),\quad
\im \omega_2+\ordo(1)\le|\omega'_{2,N}|\le 2|\omega_1|+|\omega_1|,
\]
as $N\to+\infty$, and since the area is of the fundamental rhombus is fixed,
it follows that the opening angle is controlled (bounded away from $0$ from
below as well as from above by $\pi$). It {is clear that} whenever
the geometry of the torus is controlled, as {it is} in our case,
{then there is a corresponding uniform} bound of the Green function from
above. For instance, we can simply express the Green function in terms of
the corresponding Weierstrass sigma function to obtain the bound.
\end{proof}

\begin{thm}\label{th:atom:torus1}
With the $N$ points $\zeta_1,\cdots,\zeta_N$ representing the set $\Lambda_N/\Lambda$
on the torus $\C/\Lambda$, we have that
\begin{equation*}
1\le \frac{2}{\pi}
\int_{\C/\Lambda} \exp \bigg(\sum_{j=1}^{N} \beta U_{\C/\Lambda}(\cdot,\zeta_j)\bigg)\dA
\le K^\beta,
\end{equation*}
holds for all $0\le\beta<+\infty$  for some constant $K>1$ which only depends on
the lattice $\Lambda$.
\end{thm}

\begin{proof}
In view of Lemma \ref{lem_01}, we have that
\begin{multline*}
\int_{\C/\Lambda} \exp \bigg(\sum_{j=1}^{N} \beta U_{\C/\Lambda}(\cdot,\zeta_j)\bigg)\dA
=
\int_{\C/\Lambda} \exp \left( \beta U_{\C/\Lambda_{N}}(\cdot,0)\right)\dA
\\
=N\int_{\C/\Lambda_N} \exp \left( \beta U_{\C/\Lambda_{N}}(\cdot,0)\right)\dA
=\int_{\C/\Lambda'_N} \exp \left( \beta U_{\C/\Lambda'_{N}}(\cdot,0)\right)\dA,
\end{multline*}
where $\Lambda'_N$ is the renormalized lattice of Lemma \ref{lem:renormalized1}.
In the last step, we used the scale-invariance property that 
\[
U_{\C/\Lambda_{N}}(z,0) = U_{\C/\Lambda'_N}(z\sqrt{N},0),\qquad z\in\C.
\]
The estimate of the theorem now follows from Lemma \ref{lem:renormalized1} with
the constant $K=\e^M$.
\end{proof}

\begin{proof}[Proof of Theorem \ref{mt:2}]

The assertion is an immediate corollary to Theorem \ref{th:atom:torus1}.
\end{proof}



\section{Toroidal pseudopolynomials}\label{sec_TP}

\subsection{Periodicity properties of toroidal pseudopolynomials}

We follow the notation in Subsection \ref{ss:pseudopols} and introduce
the monic first degree toroidal pseudopolynomials 
\begin{multline}
\pi_\alpha(z)=\e^{-|z|^2}\Tope_{-\alpha}(\sigma_\ast)(z)=\e^{\imag2\im \bar\alpha z}
\,\e^{-|z-\alpha|^2}\sigma_\ast(z-\alpha)
\\
=\e^{\imag2\im \bar\alpha z}\pi_0(z-\alpha).
\label{eq:monic1stdegree}
\end{multline}
As such, they are not functions on the torus $\C/\Lambda$, but rather instead
\emph{sections}. For $j=1,2$, they enjoy the periodicity-type properties
\begin{equation}
\label{eq:torpol1}
\pi_{\alpha+\omega_j}(z)=-\,\e^{\imag 2\im\bar\omega_j \alpha}\pi_\alpha(z)
\end{equation}
and
\begin{equation}
\label{eq:torpol2}
\pi_\alpha(z+\omega_j)=\chi_j(\alpha,z)\pi_\alpha(z),
\end{equation}
where $\chi_j$ stands for the \emph{phase factor} 
\begin{equation}
\chi_j(\alpha,z):=
-\e^{\imag 2\im[\bar\omega_j(z-2\alpha)]}.
\label{eq:torpol2.1}
\end{equation}
In particular, \eqref{eq:torpol1} expresses that $\pi_{\alpha+\omega_j}$ is 
an unimodular constant multiple of $\pi_\alpha$. 

\begin{defn}
The \emph{monic pseudo\-polynomials of degree $N$} on the torus $\C/\Lambda$ 
are the functions 
\[
\Pi_{\balpha}(z):=\pi_{\alpha_1}(z)\cdots\pi_{\alpha_N}(z),\qquad
\balpha=(\alpha_1,\ldots,\alpha_N),
\]
where the vector $\balpha=(\alpha_1,\ldots,\alpha_N)\in\C^N$ expresses
the \emph{roots} of $\Pi_{\balpha}$. More generally, the
\emph{pseudopolynomials of degree $N$} are the functions $\lambda\Pi_{\balpha}$,
where $\lambda\in\C$ is a constant.  
\end{defn}

The higher degree monic pseudopolynomials inherit periodicity-type properties
from the monic first degree pseudopolynomials. So, from \eqref{eq:torpol1},
we obtain
\begin{equation}\label{eq:torpol1.1}
\Pi_{\balpha+\balpha'}(z)=
\mathrm{sgn}(\balpha')\,  
\e^{\imag2\im\langle\balpha,\balpha'\rangle}\Pi_{\balpha}(z),
\qquad \balpha'\in\Lambda^N,
\end{equation}
where we use the vector notation $\balpha=(\alpha_1,\ldots,\alpha_N)$ and
$\balpha'=(\alpha'_1,\ldots,\alpha'_N)$ and associated inner product
\[
\langle\balpha,\balpha'\rangle:=\alpha_1\bar\alpha'_1+\cdots+\alpha_N
\bar\alpha'_N.
\]
Here, $\mathrm{sgn}(\balpha')$ takes values in the set $\{1,-1\}$, according to
the multiplicative rules
\[
\mathrm{sgn}(\balpha')=\mathrm{sgn}(\alpha'_1)\cdots\mathrm{sgn}(\alpha'_N)
\]
where the initial input is that for $k_1,k_2\in\Z$,
\[
\mathrm{sgn}(k_1\omega_1+k_2\omega_2)=(-1)^{k_1+k_2},\qquad j=1,2.
\]  
On the other hand, we obtain from \eqref{eq:torpol2} that for $j=1,2$,
\begin{equation}
\Pi_{\balpha}(z+\omega_j)=
\chi_j(\balpha,z)\Pi_{\balpha}(z),
\label{eq:phasefactor5}
\end{equation}
with the multiindex phase factor
\begin{equation}
\chi_j(\balpha,z):=(-1)^N
\e^{\imag 2\im[\bar\omega_j(Nz-2\alpha_1-\cdots-2\alpha_N)]}.
\label{eq:phasefactor5.1}
\end{equation}
It follows from \cref{eq:torpol1.1} that apart from an unimodular constant
multiple, $\Pi_{\balpha}$ is $\Lambda^N$-periodic with respect to 
the variable $\balpha\in\C^N$. Moreover, the modulus $|\Pi_{\balpha}(z)|$ is
$\Lambda$-periodic with respect to $z$ and $\Lambda^N$-periodic with respect to
$\balpha$.


\subsection{Linear structure: fibers of toroidal pseudopolynomials}
We consider the manifold of \emph{toroidal pseudopolynomials} of degree 
$N$, denoted by $\psi\mathrm{pol}_N(\C/\Lambda)$, which consists of all 
functions $\lambda\Pi_{\balpha}$ where $\lambda\in\C$ and 
$\balpha=(\alpha_1,\ldots,\alpha_N)\in\C^N$. Clearly, it has $\C$-dimension $N+1$.
As defined, all we can do is to move zero sets $\balpha$ and the ``leading coefficient''
parameter $\lambda$. However, sometimes the linear combination of two toroidal
pseudopolynomials is again a toroidal pseudopolynomial. To describe when this happens,
we introduce \emph{fibers} in $\psi\mathrm{pol}_N(\C/\Lambda)$.
Given a point $\gamma\in\C$, we consider the fiber
\[
\psi\mathrm{pol}_N^{\gamma}(\C/\Lambda):=\big\{\lambda\Pi_{\balpha}:
\,\alpha_1+\cdots+\alpha_N-\gamma\in\Lambda,\,\,\lambda\in\C\big\}.
\]
In view of the defining condition, we can add an element of the lattice $\Lambda$
to the parameter $\gamma$ without altering the fiber. This means that we may think of
$\gamma$ as a point on the torus: $\gamma\in\C/\Lambda$.
Since in the definition of the fiber, we reduce the manifold by one (complex) condition,
the $\C$-dimension of the fiber drops down to $N$.

It turns out that it is possible to completely characterize
functions in the fiber in terms of two simple conditions: (1) a differential
equation, and (2) periodicity-type conditions involving the two generators of the lattice
$\Lambda$. In particular, the fiber is a $\C$-linear manifold of dimension $N$.

\begin{prop}
A function $F:\C\to\C$ belongs to the fiber $\psi\mathrm{pol}_N^{\gamma}(\C/\Lambda)$
if and only if it is $C^1$-smooth and solves the differential equation
\begin{equation} 
\bar\partial_z F+NzF=0
\tag{a}
\end{equation}
while it enjoys the periodicity-type properties
\begin{equation}
F(z+\omega_j)=(-1)^N\e^{\imag 2\im[\bar\omega_j(Nz-2\gamma)]}F(z),\qquad j=1,2. 
\tag{b}
\end{equation}
\label{prop:fiber}
\end{prop}

\begin{rem}
The phase factor $\chi_j(\gamma,N,z)=\e^{\imag 2\im[\bar\omega_j(Nz-2\gamma)]}$
appearing in condition (b) only depends on $\gamma\in\C/\Lambda$, since
$\chi_j(\gamma+\omega_k,N,z)=\chi_j(\gamma,N,z)$ holds for $j=1,2$ and $k=1,2$,
as a consequence of our normalizations $\omega_1>0$ and the area condition
\eqref{eq:areacond1}. 
\end{rem}

\begin{proof}[Proof of Proposition \ref{prop:fiber}]
It is a consequence of \eqref{eq:monic1stdegree} that
\[
\Pi_{\balpha}(z)=\e^{-N|z|^2}\prod_{j=1}^{N}\Tope_{-\alpha_j}(\sigma_\ast)(z)
\]
and hence that
\[
\bar\partial_z\Pi_{\balpha}+Nz\Pi_{\balpha}=0
\]
holds. Consequently, any $F=\lambda\Pi_{\balpha}$ with 
constant $\lambda\in\C$ has the property (a). It is also automatically $C^1$-smooth.
Moreover, it follows from the
periodicity-type property \eqref{eq:phasefactor5} together with \eqref{eq:phasefactor5.1}
and the definition of the fiber $\psi\mathrm{pol}_N^{\gamma}(\C/\Lambda)$ that property
(b) holds for $F=\lambda\Pi_{\balpha}$ as well provided that
$\lambda\Pi_{\balpha}\in\psi\mathrm{pol}_N^{\gamma}(\C/\Lambda)$.

On the other hand, suppose instead that $F$ is $C^1$-smooth such that the properties
(a) and (b) hold. Then in view of (a),
\[
F(z)=\e^{-N|z|^2}G(z),
\]
where $G$ is entire. Next, choose a vector $\balpha\in\C^N$ such that
\[
\alpha_1+\cdots+\alpha_N-\gamma\in\Lambda,
\]
so that
$\Pi_{\balpha}\in\psi\mathrm{pol}_N^{\gamma}(\C/\Lambda)$. Then by condition (a), the ratio
$R:=F/\Pi_{\balpha}$ is a meromorphic function in $\C$, an in view of (b) it is also
$\Lambda$-periodic: $R(z+\omega_j)=R(z)$ for $j=1,2$. In other words, $R$ is a meromorphic
function on the torus $\C/\Lambda$ with poles at the points $\alpha_1,\ldots,\alpha_N$.
Such meromorphic functions have exactly as many zeros as they have poles (see, for
instance, Proposition 2.5.7 of \cite{NapierRamachandran2011}), and since
$F=R\Pi_{\balpha}$ the zeros of $R$ correspond to zeros of $F$.
It is now permitted to have the vector of roots $\balpha$ from $\C/\Lambda$ such that
we get complete cancellation of roots and zeros in the meromorphic function $R$, in which
case $R$ is constant, say $R=\lambda$. For that choice of $\balpha$ then, $F=\lambda
\Pi_{\balpha}\in\psi\mathrm{pol}_N^{\gamma}(\C/\Lambda)$, as claimed. 
\end{proof}

\subsection{The manifold of toroidal psudopolynomials}
While the separate fibers $\psi\mathrm{pol}_N^\gamma(\C/\Lambda)$ constitute
$N$-dimensional linear manifolds for each $\gamma\in\C/\Lambda$, the entire manifold
of toroidal pseudopolynomials $\psi\mathrm{pol}_N(\C/\Lambda)$ is nonlinear, as
the phase factors in the periodicity-type conditions depend  essentially
on the parameter $\gamma$. The fibers allow us to model the manifold
$\psi\mathrm{pol}_N(\C/\Lambda)$ as an $N$-dimensional vector bundle
(the fiber $\psi\mathrm{pol}_N^\gamma(\C/\Lambda)$) over each point on the torus
$\gamma\in\C/\Lambda$. Let us see where this leads us in terms of the definition
of $\psi\mathrm{pol}_N(\C/\Lambda)$ as consisting of $\lambda\Pi_{\balpha}$, where
$\lambda\in\C$ and $\balpha\in\C^N$. The information content of $\lambda\Pi_{\balpha}$
can be reduced down to the tuple $(\balpha,\lambda)\in\C^{N+1}$. We should describe
an equivalence relation $\cong$ among such tuples that give rise to the same element
$\lambda\Pi_{\balpha}$. A first observation is that since $\lambda$
acts multiplicatively, we have $(\balpha,0)\cong(\balpha',0)$ for all
$\balpha,\balpha'\in\C^N$. Moreover, the periodicity-type property \eqref{eq:torpol1.1}
amounts to the identification
\[
(\balpha,\lambda)\cong (\balpha+\balpha',\mathrm{sgn}(\balpha')\,
\e^{-2\imag\im\langle\balpha,\balpha'\rangle}\lambda),\qquad \balpha'\in\Lambda^N.
\]
Note that this induces a unimodular twist in the second coordinate $\lambda$ as
$\balpha$ moves in the $N$-torus $\C^N/\Lambda^N$.
In the parameter space, each fiber $\psi\mathrm{pol}_N^\gamma(\C/\Lambda)$ of course
corresponds to the set 
\[
\big\{(\balpha,\lambda):\,\alpha_1+\cdots+\alpha_N-\gamma\in\Lambda\big\}.
\]

\subsection{Norm estimate of toroidal pseudopolynomials}
We should supply the proof of Proposition \ref{mt:3}.

\begin{proof}[Proof of Proposition \ref{mt:3}]
Given the relation \eqref{eq:Lbetanorm2.01}, the assertion of the proposition just
restates the estimate \eqref{eq:04} in a slightly different form.
\end{proof}

\end{document}